\RequirePackage{ifpdf}
\ifpdf 
\documentclass[pdftex]{sigma}
\else
\documentclass{sigma}
\fi

\begin{document}
\allowdisplaybreaks

\renewcommand{\PaperNumber}{003}

\FirstPageHeading

\ShortArticleName{Measures on Current Groups}

\ArticleName{Heat Kernel Measure on Central Extension\\ of Current
Groups in any Dimension}

\Author{R\'emi L\'EANDRE}

\AuthorNameForHeading{R. L\'eandre}

\Address{Institut de Math\'ematiques de Bourgogne, Universit\'e de
Bourgogne, 21000 Dijon,  France}

\Email{\href{mailto:Remi.leandre@u-bourgogne.fr}{Remi.leandre@u-bourgogne.fr}}

\ArticleDates{Received October 30, 2005, in final form January 13, 2006; Published online January 13, 2006}

\Abstract{We define measures on central extension of current
groups in any dimension by using infinite dimensional Brownian motion.}


\Classification{22E65; 60G60}

\section{Introduction}
If we consider a smooth loop group, the basical central extension
associated to a suitable Kac--Moody cocycle plays a big role
 in mathematical physics \cite{3,11,21,24}. L\'eandre
 has defined the space of $L^2$ functionals on a continuous Kac--Moody group,
 by using the Brownian bridge measure on the basis~\cite{16}
 and deduced the so-called energy representation of the smooth Kac--Moody group on it.
 This extends the very well known representation of a
 loop group of Albeverio--Hoegh--Krohn~\cite{2}.

Etingof--Frenkel~\cite{13} and Frenkel--Khesin~\cite{14} extend
these considerations to the case where the parameter space is two dimensional.
They consider a
 compact Riemannian surface $\Sigma$ and consider the set of
  smooth maps from $\Sigma$ into a compact simply connected
  Lie group $G$. We call $C_r(\Sigma;G)$ the
  space of $C^r$ maps from $\Sigma$ into $G$
  and $C_\infty(\Sigma;G)$ the space of smooth maps
  from $\Sigma$ into $G$. They consider the universal cover
   $\tilde{C}_\infty(\Sigma;G)$ of it and construct
   a central extension by the Jacobian~$J$ of $\Sigma$ of it $\hat{C}_\infty(\Sigma;G)$ (see \cite{7, 8, 25}
 for related works).

We can repeat this construction if $r>s$ big enough for
$C_r(\Sigma;G)$. We get the universal cover
$\tilde{C}_r(\Sigma;G)$
 and the central extension by the Jacobian $J$ of $\Sigma$
 of $\tilde{C}_r(\Sigma;G)$ denoted by $\hat{C}_r(\Sigma;G)$.

By using Airault--Malliavin construction of the Brownian motion on
a loop group \cite{1,9}, we have  defined in \cite{19} a
probability measure on
 $\tilde{C}_r(\Sigma;J)$, and since the Jacobian is compact,
  we can define in \cite{19}
 a probability measure on $\hat{C}_r(\Sigma;G)$.

Maier--Neeb \cite{20} have defined the universal central extension
of a current group $C_\infty(M;G)$ where $M$ is any compact
manifold. The extension is done by a quotient of a certain space
of differential form on $M$ by a lattice.
 We remark that the Maier--Neeb procedure can be used if we replace this infinite dimensional space of forms by
 the de Rham cohomology groups $H(M:{\rm Lie}\,G)$ of $M$ with values in ${\rm Lie}\,G$.
 Doing this,  we get a central extension by a finite dimensional Abelian groups instead of an infinite dimensional
Abelian group. On the current group $C_r(M;G)$ of $C^r$ maps from
$M$ into the considered compact connected Lie group $G$,
 we use heat-kernel measure deduced from the Airault--Malliavin equation,
 and since we get a central extension $\hat{C}_r(M;G)$ by a~finite dimensional group $Z$,
we get a measure on the central extension of the current group.
Let us recall that studies of the Brownian motion on infinite
dimensional manifold have a long history (see works of Kuo \cite{15},
Belopolskaya--Daletskii \cite{6,12}, Baxendale \cite{4,5}, etc.).

Let us remark that this procedure of getting a random field by
adding extra-time is very classical in theoretical physics,
 in the so called programme of stochastic-quantization of Parisi--Wu~\cite{23},
 which uses an infinite-dimensional Langevin equation.
 Instead to use here the Langevin equation, we use  the more tractable Airault--Malliavin equation,
 that represents infinite dimensional Brownian motion on a current group.

\section{A measure on the current group in any dimension}

We consider $C_r(M;G)$ endowed with its $C^r$ topology. The
parameter space $M$ is supposed compact and the Lie group $G$ is
supposed compact, simple and simply connected.
  We consider the set of continuous paths from $[0,1]$ into $C_r(M;G)$
  $t \rightarrow g_t(\cdot)$, where $S\in M \rightarrow g_t(S)$ belongs to $C_r(M
;G)$ and $g_0(S) = e$. We denote $P(C_r(M;G))$ this path space.

Let us consider the Hilbert space $H$ of maps $h$ from $M$ into
${\rm Lie}\, G$ defined as follows:
\begin{gather*}
\int_\Sigma\langle (\Delta^k+1)h, h\rangle dS = \Vert h \Vert_H^2,
\end{gather*}
where $\Delta$ is the Laplace Beltrami operator on $M$ and $dS$
the Riemannian element on $M$ endowed with a Riemannian structure.

We consider the Brownian motion $B_t(\cdot)$ with values in $H$.

We consider the Airault--Malliavin equation (in Stratonovitch sense):
\begin{gather*}
dg_t(S) = g_t(S)dB_t(S),\qquad g_0(S) = e.
\end{gather*}
Let us recall (see \cite{17}):

\begin{theorem} If $k$ is  enough big, $t \rightarrow\{S \rightarrow g_t(S)\}$
defines a random element of $P(C_r(M;G))$.
\end{theorem}

We denote by $\mu$ the heat-kernel measure  $C_r(M;G)$: it is the
law of the $C^r$ random field $S \rightarrow g_1(S)$.
 It is in fact a probability law on the connected component of the identity
 $C_r(M;G)_e$ in the current group.

\section[A brief review of Maier-Neeb theory]{A brief review of Maier--Neeb theory}
Let us consider $\Pi_2(C_r(M;G)_e)$ the second fundamental group
of the identity in the current group for $r>1$. The Lie algebra of
this current group is $C_r(M; {\rm Lie}\, G)$ the space of $C^r$
maps from $M$ into the Lie algebra ${\rm Lie}\, G$ of $G$
\cite{22}. We introduce the canonical Killing form $k$ on ${\rm
Lie}\, G$.

 $\Omega^i(M; {\rm Lie}\, G)$ denotes the space of $C^{r-1}$
 forms of degree $i$ on $M$ with values in ${\rm Lie}\, G$. Follo\-wing~\cite{20}, we
introduce the left-invariant 2-form $\Omega$ on $C_r(M;G)$ with
values in the space of forms $Y = \Omega^1(M;{\rm Lie}\,
G)/d\Omega^0(M; {\rm Lie}\, G)$ which associates
\begin{gather*}
k(\eta, d\eta_1).
\end{gather*}
to $(\eta, \eta_1)$, elements of the Lie algebra of the current
group.

For that, let us recall that the Lie algebra of the current group
is the set of $C^r$ maps $\eta$ from the manifold into the Lie
algebra of $G$. $d\eta$ is a $C^{r-1}$ 1-form into the Lie algebra
of $G$. Therefore $k(\eta,d\eta_1)$ appears as a $C^{r-1}$ 1-form
with values in the Lie algebra of $G$.
 Moreover
\begin{gather*}
dk(\eta, \eta_1) = k(d\eta, \eta_1) + k(\eta, d\eta_1).
\end{gather*}
This explains the introduction of the quotient in $Y$. Following
the terminology of~\cite{20}, we consider the period map $P_1$
which to $\sigma$ belonging to
$\Pi_2(C_r(M;G)_e)$ associates $\int_\sigma \Omega$. Apparently~$P_1$ takes its values in $Y$, but
in fact, the period map takes its values in a lattice $L$ of
$H^1(M; {\rm Lie}\, G)$. It is defined on
 $\Pi_2(C_r(M;G)_e)$ since $\Omega$ is closed
for the de Rham differential on the current group, as it is
left-invariant and closed and it is a 2-cocycle in the Lie algebra of
the current group~\cite{20}. We consider the Abelian group $Z =
H^1(M; {\rm Lie}\, G)/L$. $Z$ is {\it of finite dimension}.

We would like to apply Theorem III.5 of~\cite{20}. We remark that
the map $P_2$ considered as taking its values in $Y/L$ is still
equal to 0 when it is considered by taking its values in $H^1(M;
{\rm Lie}\, G)/L$.

 We deduce the following theorem:

\begin{theorem} We get a central extension $\hat{C}_r(M;G)$ by $Z$ of the
current group $C_r(M;G)_e$ if $r>1$.
\end{theorem}

Since $Z$ is of finite dimension, we can consider the Haar measure
on $Z$. We deduce from $\mu$ a~measure $\hat{\mu}$ on
$\hat{C}_r(M;G)$.

\begin{remark} Instead of considering $C_r(M; {\rm Lie}\, G)$,
we can consider $W_{\theta,p}(M;{\rm Lie}\, G)$,
 some convenient Sobolev--Slobodetsky spaces of maps from $M$
 into ${\rm Lie}\, G$. We can deduce a central extension $\hat{C}_{\theta, p}(M;G)$ of
the Sobolev--Slobodetsky current group $C_{\theta,p}(M;G)_e$. This
will give us an example of Brzezniak--Elworthy theory, which works
for the construction of diffusion processes on
infinite-dimensional manifolds modelled on M-2 Banach spaces,
since Sobolev--Slobodesty spaces are M-2 Banach spaces
\cite{9,10,18}. We consider a Brownian motion $B^1_t$ with values
in the finite dimensional Lie algebra of $Z$ and $ \hat{B}_t =
(B_t(\cdot), B^1_t)$ where $B_t(\cdot)$ is the Brownian motion
in~$H$ considered in the Section~2. Then, following the ideas of
Brzezniak--Elworthy,
 we can consider the stochastic differential equation on $\hat{C}_{\theta,p}(M;G)$
 (in Stratonovitch sense):
\begin{gather*}
d\hat{g}_t(\cdot) = \hat{g}_t(\cdot)d\hat{B}_t.
\end{gather*}
\end{remark}

\subsection*{Acknowledgements}
The Author  thanks Professor K.H.~Neeb for helpful discussions.

\LastPageEnding
\end{document}